# Sustainability and coordination in a socially responsible supply chain using a combined incentive contract and a social marketing strategy

Mahdi Ebrahimzadeh-Afrouzi, Masoud Asadpour Ahmadchali

School of Industrial Engineering, Iran University of Science and Technology, Tehran, Iran

[mehdi.eh90@gmail.com](mailto:mehdi.eh90@gmail.com), [ma.asadpour1990@gmail.com](mailto:ma.asadpour1990@gmail.com)

## Abstract

In recent decades, due to the growing concerns for sustainable development, supply chains seek to invest in social sustainability issues to seize more market share in today's competitive business environment. This study aims to develop a coordination scheme for a manufacturer-retailer supply chain (SC) contributing in social donation (SD) activity under a cause-related marketing (CRM) campaign. In the presence of consumer social awareness (CSA), the manufacturer notices consumers through some activities (i.e. labelling) that he participates in a CRM campaign by donating a proportion of the retail price to a cause whenever a consumer makes a purchase. In this study, the market demand is dependent on the retail price, the retailer's stock level and donation size. The proposed problem is designed under three decision-making systems. Firstly, a decentralized decision-making system (traditional structure), where the SC's members aim to optimize their own profits regardless of the other member's profitability, is investigated. Then, the problem is designed under a centralized decision-making system to obtain the best values of retail price and replenishment decisions from the entire SC perspective. Afterwards, an incentive mechanism based on a cost and revenue-sharing (RCS) factor is developed in the coordination system to persuade the SC members to accept the optimal results of the centralized system without suffering any profit loss. Moreover, the surplus profit obtained in the centralized system is divided between the members based on their bargaining power. The numerical investigations and the blocked decision-making on SD activity are presented to evaluate the proposed model. Not only does the proposed coordination model increase the SC members' profit, but it is also desirable in achieving a more socially responsible SC.

**Keywords:** Sustainable supply chain coordination, Pricing, Inventory-dependent demand, Revenue- and cost-sharing contract, cause-related marketing

## 1. Introduction

In recent decades, the growing development of the industry and the intense competition of companies to increase their profits have brought about adverse and destructive impacts on the environment and societies (Nigam, 2014). As a result of such detrimental effects on the world, the concepts of sustainability and sustainable development system have emerged. In this regard, customers expect companies and industrial institutes to be responsible for such consequences and consider corporate social responsibility (CSR) in their business operations. The significance of sustainable development is tied to the concept of CSR and there is a profound relation between these concepts (Kolk 2016). As consumer social awareness (CSA) increases, supply chains are required to be more socially responsible in order to survive in today's tight competition. Nowadays, governments and competitive corporations aim to persuade the downstream retailers to sell and distribute social products (Klassen 2012). A social product is something that benefits the largest number of people in the largest possible way, such as clean air, clean water, healthcare, altruism and literacy.

In line with increasing pressure on companies to be more socially responsible, cause-related marketing (CRM) has been frequently used as a popular marketing process (Galán Ladero et al 2013). Indeed, CRM is a certain kind of CSR, introduced as marketing activities by Lii and Lee, 2012, in such a way that companies donate a specific amount to a certain cause. In other words, in CRM, companies commit to donating a proportion of the sales revenue to a cause whenever a customer makes a purchase, which leads to more profitability for companies (He et al 2016). In this context, most researches substantiated that CRM can increase sales volume (Andrews et al., 2014), sales price (Leszczyc and Rothkopf, 2010), and even strengthen the relationship between companies towards CSR (Bhattacharya et al., 2009). CRM was first implemented in 1983 by American Express in a way that it organized a campaign to rebuild the Statue of Liberty and donated an amount each time to customers using their American Express cards. Later, many other companies and organizations participated in CRM campaigns so that in 2017, 2.05 million dollars were devoted to CRM in America . A clear example of these CRM campaigns are Starbucks that donated one dollar from each sale of Africa Blend coffee to support AIDS patients in Africa (Müller et al., 2014).

Donation size, which is considered as a CRM campaign design factor, can be an important factor in the success of a supply chain (SC) in the CSR area. Müller et al. (2014) expressed that donation size can have a positive impact on the choice of the product by customers, resulting in an increase in the demand for that product. Therefore, since the donation size affects the performance of SCs, controlling the proper decisions, i.e. pricing and replenishment, is of great importance in gaining profit from participation in the CRM campaign. Since the donated volume is determined by the SC's upstream member and its cost is paid by him, operating under the traditional decision-making system (decentralized system) may not lead to profitability for the whole SC or one of its members. In this vein, the coordination of the decisions is necessary to ensure the maximization of the benefit of the whole SC and its members.

Coordination mechanisms are a guaranteed effective way to avoid wasting the potential benefits that a whole SC and its members can earn. Coordination is also necessary for some industries that deal with perishable products, such as fruit, medicine, dairy, meat stores, etc. (Cachon 2003). In this article, due to the perishable nature of the products, the demand that the chain faces is considered dependent on the retailer's inventory level. Mainly, such a demand pattern is for perishable products, whose instantaneous demand has a close relationship with the instantaneous inventory of the product (Yang et al 2014). Moreover, because environmental and social costs are usually handled by manufacturers alone, choosing a proper coordination mechanism is very important, especially for manufacturers.

In this article, a two-echelon SC, including one manufacturer and one retailer, encountering socially aware customers is investigated. In order to choose its brand by socially aware customers, the manufacturer participates in the CRM campaign and notices customers through activities such as labelling that he donates a percentage of the sales price for social activities. In the proposed model, the market demand depends on retail price, donation size and the retailer's inventory level. As a result of the impact of donation size on the SC, simultaneous integrated decision making on pricing and replenishment policies are investigated. Firstly, the problem is modelled under the decentralized and centralized decision-making systems. Under the decentralized decision-making system, in which the decisions are made independently regardless of the desirability of the other member, the retailer aims to optimize the pricing and replenishment policies while the manufacturer decides on the number of inventory replenishment. Clearly, the centralized system as an integrated channel seeks to optimize the profitability of the whole SC regardless of the profit of SC members; therefore, decisions made under the centralized system are the optimal solutions (maximum profit) for the whole SC. In this regard, although investing in CSR campaign increases the market demand and the whole SC's profit, this investment is costly for the manufacturer. Therefore, it may cause a financial loss to a member while leading to much more profit for the whole SC by making decision under the centralized system. Under such a circumstance, it is necessary to select an incentive mechanism under a coordination system to coordinate the decisions of the whole SC and its members in order to reach the maximum profit obtained from operating under the centralized system so that none of the members suffers losses by participating in the CSR campaign. In this paper, in order to motivate the SC members to coordinate their decisions towards the optimal solutions of the centralized system, the revenue and cost-sharing (RCS) mechanism as a combined contract is proposed. Despites, a profit-sharing strategy modelled based on the SC members' bargaining power is applied to fairly allocate the surplus profit gained from the centralized system.

The contribution of this study is on the social SC coordination participating in CRM campaigns by designing a combined revenue and cost-sharing contract for simultaneous coordination of the pricing, donation size and replenishment policies. The main purpose of this article is to answer the following research questions.

- How three decision-making systems (i.e. decentralized, centralized and coordinated) can affect the optimal conditions of the whole SC and its members when the demand is dependent on retail price, donation size and the retailer's inventory level?
- What is the optimal value of retail price, order quantity and the number of manufacturer's replenishment under decentralized and centralized decision-making systems?
- Is it feasible to implement a coordinating incentive scheme in order to encourage the whole SC and its members to shift their decisions in the decentralized system to the centralized system's optimal solutions?
- To what extent does participation in the CRM campaign affect the performance and profitability of SC and its members?

In the following, the related literature review is provided in Section 2. Problem definition and computational models for different types of decision-making systems are represented in Section 3. Then, we provide a set of illustrative examples and sensitivity analyses in Section 4. To analyze the performance of the SC in the CRM campaigns, a blocked decision-making structure on social donation size with related numerical results are presented in Section 5. Moreover, the managerial insights are presented in Section 6. Finally, Section 7 describes the conclusions of this study and potential extensions for the model in future studies.

## 2. Theoretical background

This section intentionally left blank

## 3. Problem definition and mathematical models

This paper investigates a manufacturer-retailer supply chain (SC), in which the manufacturer aims to participate in social donation (SD) activity under a cause-related marketing (CRM) strategy to increase the sales volume thanks to the customer social awareness (CSA). In other words, the manufacturer decides to donate a ration $\theta$ of the retail price per each sold item to an SD activity and obviously notices $\theta$ to customers by sticking some labels on their products. The CRM costs are all handled by the manufacturer; therefore, no increase occurs in the amount of wholesale price $v$. The manufacturer uses a production-replenishment system for loading inventories and decides on the number of replenishments $n$ as a decision variable to determine the quantity of production. The retailer follows the EOQ system and makes decision on retail price $p$ and order quantity $Q$ to replenish its inventory level.

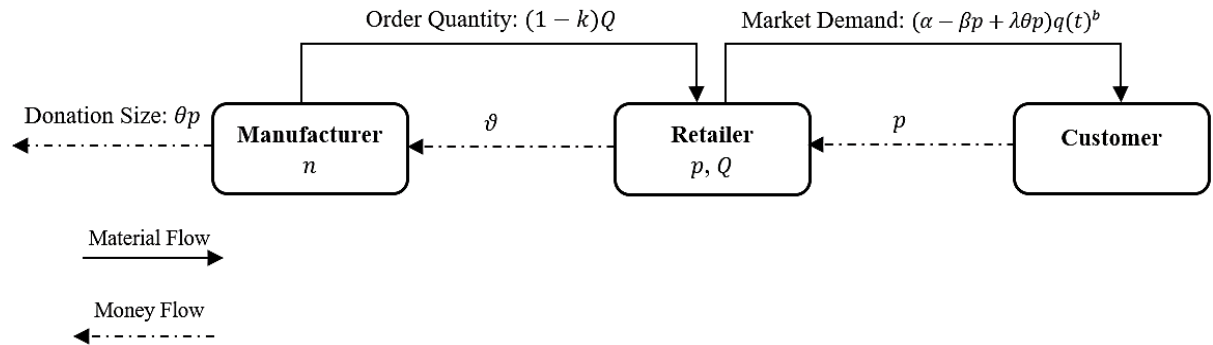

**Fig 1**. SC investigated in this paper.

Table 1. The notations in this paper.

| | Description |
|---|---|
| • **Decision variables** | |
| $Q$ | Retailer's order quantity |
| $p$ | Retail price per unit |
| $n$ | Number of delivery from manufacturer to the retailer per setup run, a positive integer |
| $\mu$ | The revenue fraction of the retailer that generated in the incentive contract, $0 < \mu < 1$ |
| $v_{co}$ | Discounted wholesale price in the incentive contract |
| • **Parameters** | |
| $q(t)$ | Retailer's inventory level at time period $t$ |
| $T_r$ | Replenishment cycle length of the retailer |
| $T$ | Replenishment cycle length |
| $\vartheta$ | Manufacturer's wholesale price per unit |
| $m$ | Production cost of the manufacturer per unit, $(m < \vartheta < p)$ |
| $\theta$ | Portion of retail price that the manufacturer donates in SD activity, $0 \leq \theta < 1$ |
| $k$ | Retailer's reorder point, $0 < k < 1$ |
| $R$ | Manufacturer's production rate |
| $A_r$ | Retailer's ordering cost |
| $A_m$ | Manufacturer's setup cost |
| $h_r$ | Retailer's holding cost per unit for per unit time |
| $h_m$ | Manufacturer's holding cost per unit for per unit time |
| $\xi$ | Bargaining power of the retailer versus the supplier, $0 < \xi < 1$ |
| $\Pi_r$ | Retailer's average profit |
| $\Pi_m$ | Manufacturer's average profit |
| $\Pi_{sc}$ | Supply chain average profit |

Note: Subscripts d, c and co imply the decentralized, centralized and coordinated systems, respectively.

The notations are described in Table 1 and the following assumptions are considered throughout the whole paper.

### 3.1. Assumptions

- The supply chain consists of a single-manufacturer and a single-retailer trading a single type of product.
- The market demand rate is dependent on retail price $p$, donation size $\theta p$ and the retailer's instantaneous inventory level $q(t)$, i.e., $D(p,t) = (\alpha - \beta p + \lambda \theta p) q(t)^b$, where $\alpha > 0$ is a market potential scale, $\beta > 0$ is the retail price-sensitive parameter of demand, $\lambda > 0$ is customer social awareness (CSA) level, and $0 < b < 1$ is the demand elasticity with respect to the retailer's current inventory level. In the first part of the demand function, according to Müller et al. (2014), there is a strong relation between the market demand and donation size, meaning that by increasing the donation size announced to customers by labelling on the product, the market demand increases much more quickly. Such a type of demand pattern can be seen commonly in various SC models particularly in social and green SCs in the literature.
- The retailer's inventory is continuously reviewed at an interval of time $T_r$. An order with a size of $(1-k)Q$ is replenished whenever the retailer's inventory level falls to the reorder point $kQ$ (see Figure 2).
- The manufacturer uses order multiplier $n$ to replenish its inventory. The retailer places an order $(1-k)Q$, and the manufacturer produces with a finite production rate $R(> D)$ in a single setup $n(1-k)Q$ but transfers the entire lot to the retailer over $n$ shipments of equal size.
- Shortages are not allowed to occur.
- The lead time is zero.

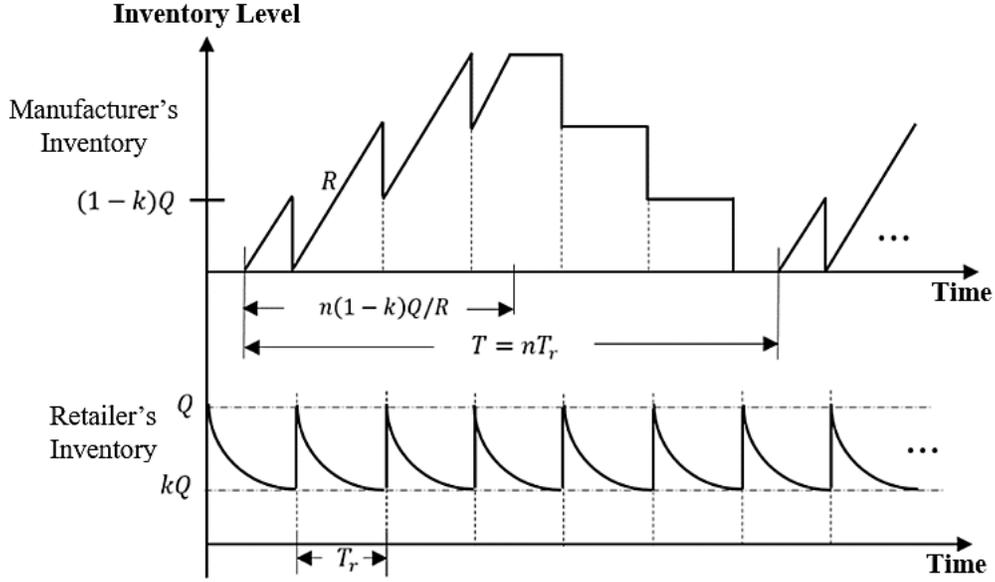

**Fig 2.** The inventory pattern for the retailer and the manufacturer.

Since the market demand rate is equal to the decrease in the inventory level, the retailer's inventory level $q(t)$ can be described by the following equation:

$$\frac{\partial q(t)}{\partial t} = -(\alpha - \beta p + \lambda\theta p)q(t)^b \qquad 0 \leq t \leq T_r \qquad (1)$$

By integrating Eq. (1), we have:

$$\int_0^t q(t)^{-b}\, dt = -\int_0^t (\alpha - \beta p + \lambda\theta p)dt \qquad (2)$$

and

$$q(t)^{1-b} - q(0)^{1-b} = -(\alpha - \beta p + \lambda\theta p)(1-b)t. \qquad (3)$$

Since $q(0) = Q$ and the retailer's inventory level decreases to reorder point at time $T_r$, $q(T_r) = kQ$, we get

$$q(t) = [Q^{1-b} - (\alpha - \beta p + \lambda\theta p)(1-b)t]^{\frac{1}{1-b}} \qquad (4)$$

and

$$T_r = \frac{(1-k^{1-b})Q^{1-b}}{(1-b)(\alpha-\beta p+\lambda\theta p)}, \qquad p < \frac{\alpha}{\beta - \lambda\theta} \qquad (5)$$

### 3.2. Decentralized system

In this section, under decentralized system as a traditional decision-making system, decisions are optimized independently by SC members based on their own interests to maximize their own profits without considering the whole SC profit.

#### 3.2.1. The retailer's model

The retail price $p$ and order quantity $Q$ are determined by the retailer to maximize its profit without considering the whole SC profit. The retailer's profit function is modeled as follows,

$$\text{Maximize } \Pi_{r/d}(p,Q) = \frac{1}{T_r}\left[(p-v)(Q-kQ) - A_r - h_r \int_0^{T_r} q(t)dt\right] = \frac{1}{T_r}\left[(p-v)(1-k)Q - A_r - \frac{(1-k^{2-b})h_r}{(\alpha-\beta p+\lambda\theta p)(2-b)}Q^{2-b}\right] \qquad p < \frac{\alpha}{\beta-\lambda\theta} \qquad (6)$$

The elements of Eq. (6) are sales revenue, ordering cost and holding cost, respectively. Substituting Eq. (5) into Eq. (6) and simplifying, we obtain

$$\Pi_{r/d}(p.Q) = \frac{(1-b)(\alpha-\beta p+\lambda\theta p)}{(1-k^{1-b})}[(p-v)(1-k)Q^b - A_r Q^{b-1}] - \frac{(1-b)(1-k^{2-b})h_r}{(2-b)(1-k^{1-b})}Q \qquad p < \frac{\alpha}{\beta-\lambda\theta} \qquad (7)$$

Since the objective function is a power function, for convenience, the optimality conditions can be determined in a manner similar to what was done by Sheen and Tsao (2007). In their method, first the optimality conditions at $p = p_0$ for a given $Q$ is attained by satisfying the first and second order derivative conditions w.r.t. $p$ (i.e. $\partial\Pi_{r/d}/\partial p = 0$, $\partial^2\Pi_{r/d}/\partial p^2 < 0$) meaning that for a fixed $Q$, the profit function $\Pi_{r/d}$ is concave on $p$ and $p_0$ (retail price). Next, considering the concavity condition, $\Pi_{r/d}(p.Q)$ can be converted into $\Pi_{r/d}(Q)$ by substituting the resulting $p^*(Q)$ into the model. Then, the optimality conditions of $\Pi_{r/d}(Q)$ can be achieved by satisfying the first and second order derivative conditions w.r.t. $Q$.

According to above explanations, taking the second derivative of $\Pi_{r/d}(p.Q)$ w.r.t. $p$ for a given $Q$,

$$\frac{\partial^2 \Pi_{r/d}}{\partial p^2} = -\frac{2(\beta-\lambda\theta)(1-b)(1-k)}{(1-k^{1-b})}Q^b < 0 \qquad (8)$$

For a given $Q$, $\Pi_{r/d}$ is concave w.r.t. $p$ if $\theta < \frac{\beta}{\lambda}$. Hence, we can get the optimal retail price $p^*$ as a function of the order quantity $Q$ at $\partial\Pi_{r/d}/\partial p = 0$:

$$p^*(Q) = \frac{1}{2}\left[\frac{\alpha}{\beta-\lambda\theta} + v + \frac{A_r}{(1-k)Q}\right] \qquad (9)$$

Substituting $p^*(Q)$ into Eq. (7), we can get the retailer's objective function, as a function of $Q$, as follows:

$$\Pi_{r/d}(Q) = \frac{(\beta-\lambda\theta)(1-b)}{4(1-k^{1-b})}\left\{\left[\left(\frac{\alpha}{\beta-\lambda\theta}-v\right)\sqrt{(1-k)Q^b}\right) - \frac{A_r}{\sqrt{(1-k)Q^{2-b}}}\right]^2 - \frac{4(1-k^{2-b})}{(\beta-\lambda\theta)(2-b)}h_r Q\right\} \qquad v < \frac{\alpha}{\beta-\lambda\theta}, 0 \le \theta < \frac{\beta}{\lambda} < 1 \qquad (10)$$

As mentioned above, only the concave part of $\Pi_{r/d}(Q)$, which has positive values, should be taken into account (Kim et al. (1995) and Sheen and Tsao (2007)) since $\Pi_{r/d}(Q)$ is a convex-concave function on $Q$. Taking the second derivative of $\Pi_{r/d}(Q)$ w.r.t. $Q$, we will have:

$$\frac{\partial^2 \Pi_{r/d}(Q)}{\partial Q^2} = -\frac{b(\beta-\lambda\theta)(1-b)^2\left(\frac{\alpha}{(\beta-\lambda\theta)}-v\right)^2(1-k)Q^{b-2}}{4(1-k^{1-b})} - \frac{2(\beta-\lambda\theta)(1-b)^2(2-b)\left(\frac{\alpha}{(\beta-\lambda\theta)}-v\right)A_r Q^{b-3}}{4(1-k^{1-b})} + \frac{(\beta-\lambda\theta)(1-b)(2-b)(3-b)A_r^2 Q^{b-4}}{4(1-k)(1-k^{1-b})} \quad (11)$$

For simplifying, by rewriting $\partial^2 \Pi_{r/d}(Q)/\partial Q^2$, we have

$$\frac{\partial^2 \Pi_{r/d}(Q)}{\partial Q^2} = \frac{(\beta-\lambda\theta)(1-b)}{4(1-k^{1-b})} Q^{b-4}(-\tau_1 Q^2 - \tau_2 Q + \tau_3) \quad (12)$$

Where $\tau_1 = b(1-b)\left(\frac{\alpha}{(\beta-\lambda\theta)}-v\right)^2(1-k)$, $\tau_2 = 2(1-b)(2-b)\left(\frac{\alpha}{(\beta-\lambda\theta)}-v\right)A_r$, $\tau_3 = (2-b)(3-b)A_r^2/(1-k)$ and $\tau_1, \tau_2, \tau_3 > 0$.

Setting Eq. (12) equal to zero and solving it, the two saddle points can be obtained as follows:

$$Q_1 = \frac{-\tau_2 + \sqrt{\tau_2^2 + 4\tau_1\tau_3}}{2\tau_1}, \quad (13)$$

$$Q_2 = \frac{-\tau_2 - \sqrt{\tau_2^2 + 4\tau_1\tau_3}}{2\tau_1}, \quad (14)$$

Therefore

$$\frac{\partial^2 \Pi_{r/d}(Q)}{\partial Q^2} \begin{cases} >0 & \text{if } Q_2 < Q < Q_1 \\ <0 & \text{otherwise} \end{cases} \quad (15)$$

From Eq. (15), the total profit function is concave when $Q > Q_1$ and $Q < Q_2$, and convex between the two saddle points. Moreover, $Q_2 < 0$, and it can be found that $Q_1$ is a positive value when $\tau_2 < \sqrt{\tau_2^2 + 4\tau_1\tau_3}$.

Thus, the economic order quantity $Q^*$ can be obtained when $\partial \Pi_{r/d}/\partial Q = 0$:

$$\frac{\partial \Pi_{r/d}(Q)}{\partial Q} = \frac{(\beta-\lambda\theta)(1-b)}{4(1-k^{1-b})}\left[b\left(\frac{\alpha}{(\beta-\lambda\theta)}-v\right)^2(1-k)Q^{b-1} + 2(1-b)\left(\frac{\alpha}{(\beta-\lambda\theta)}-v\right)A_r Q^{b-2} - \frac{(2-b)A_r^2 Q^{b-3}}{(1-k)} - \frac{4(1-k^{2-b})}{(\beta-\lambda\theta)(2-b)}h_r\right] = 0 \quad (16)$$

Substituting $Q^*$ into Eq. (9), we can get the optimal retail price $p^*$. Thereupon, the retailer's optimal profit $\Pi_{r/d}^*(p^*, Q^*)$ can be achieved.

3.2.2. The manufacturer's model

Under decentralized system, the manufacturer operates independently but has to follow the retailer's decisions. Manufacturer makes decision on the number of replenishments ($n$) based on its total expected profit. From Fig. 2, the manufacturer's average inventory is

$$I_m = \frac{1}{T}\left[\left(n(1-k)Q\left(\frac{(1-k)Q}{R} + \frac{n-1}{n}T\right) - \frac{[n(1-k)Q]^2}{2R}\right) - \left(\frac{(1-k)QT}{n}(1 + 2 + \cdots + (n-1))\right)\right] = \frac{(1-k)Q}{2}\left[(n-1)\left(1 - \frac{n(1-k)Q}{RT}\right) + \left(\frac{n(1-k)Q}{RT}\right)\right] \quad (17)$$

Hence, the manufacturer's expected profit function per unit time can be modeled as shown in Eq. (18), which is the supplier's expected net revenue minus ordering cost, holding cost and paid donation cost per each sold item;

$$\Pi_{m/d}(n) = \frac{(v-m)n(1-k)Q}{T} - \frac{A_m}{T} - \frac{1}{2}h_m\left[(n-1)\left(1 - \frac{n(1-k)Q}{RT}\right) + \left(\frac{n(1-k)Q}{RT}\right)\right](1-k)Q - \frac{\theta pn(1-k)Q}{T} \quad (18)$$

Substituting $T = nT_r$ and $T_r = (1 - k^{1-b})Q^{1-b}/(1-b)(\alpha - \beta p + \lambda\theta p)$ in to Eq. (16), we will have

$$\Pi_{m/d}(n) = \frac{(1-b)(\alpha-\beta p+\lambda\theta p)}{(1-k^{1-b})}\left[(v - m - \theta p)(1-k)Q^b - \frac{A_m}{n}Q^{b-1}\right] - \frac{h_m(1-k)Q}{2}\left[(n-1) + \frac{(1-b)(\alpha-\beta p+\lambda\theta p)(2-n)(1-k)Q^b}{R(1-k^{1-b})}\right] \quad (19)$$

Taking the second derivative of $\Pi_{m/d}$ w.r.t. $n$, we can find that $\Pi_{m/d}$ is concave in $n$:

$$\frac{\partial^2\Pi_{m/d}(n)}{\partial n^2} = -\frac{2(1-b)(\alpha-\beta p+\lambda\theta p)A_m}{(1-k^{1-b})Q^{1-b}n^3} < 0 \quad (20)$$

Thus, the value of $n$ can be calculated by setting $\partial\Pi_{m/d}/\partial n = 0$ as:

$$n^*_{Desimal} = \sqrt{\frac{2RA_m(1-b)(\alpha-\beta p^*+\lambda\theta p^*)Q^{*b}}{h_m(1-k)Q^{*2}\left(R(1-k^{1-b})-(1-b)(\alpha-\beta p^*+\lambda\theta p^*)(1-k)Q^{*b}\right)}}. \quad (21)$$

Since $n^*_{Desimal}$ is a positive value, Eq. (21) will have a logical value when $Q < [R(1 - k^{1-b})/(1-b)(\alpha - \beta p^* + \lambda\theta p^*)(1-k)]^{1/b}$. Therefore, as $n$ is constrained to integer values, the optimal value of $n$ can be calculated as follows:

$$n^* = \begin{cases} \lfloor n^*_{Desimal} \rfloor & if \quad \Pi_{m/d}(\lfloor n^*_{Desimal} \rfloor) > \Pi_{m/d}(\lceil n^*_{Desimal} \rceil) \\ \lceil n^*_{Desimal} \rceil & if \quad \Pi_{m/d}(\lfloor n^*_{Desimal} \rfloor) \leq \Pi_{m/d}(\lceil n^*_{Desimal} \rceil) \end{cases} \quad (22)$$

Then, the expected profit of the SC under decentralized system will be

$$\Pi_{sc/d} = \Pi^*_{r/d}(p^*, Q^*) + \Pi_{m/d}(n^*) \quad (23)$$

### 3.3. Centralized system

Under this system, it is supposed that a single decision maker, i.e. an integrated firm, decides simultaneously on decision variables to optimize the whole SC profit. So in this section, the decision variables are expected to globally be optimized. The whole SC's expected profit function, which is the sum of all SC members' profit functions in decentralized system, is expressed by Eq. (24).

$$\text{Maximize } \Pi_{sc/c}(p,Q,n) = \frac{(1-b)(\alpha-\beta p+\lambda\theta p)}{(1-k^{1-b})}\bigg[(p-m-\theta p)(1-k)Q^b - \left(A_r+\frac{A_m}{n}\right)Q^{b-1}\bigg] - \frac{(1-b)(1-k^{2-b})h_r}{(2-b)(1-k^{1-b})}Q - \frac{h_m(1-k)Q}{2}\bigg[(n-1) + \frac{(1-b)(\alpha-\beta p+\lambda\theta p)(2-n)(1-k)Q^b}{R(1-k^{1-b})}\bigg] \quad p < \frac{\alpha}{\beta-\lambda\theta}, p,Q,n \geq 0 \text{ and } n \text{ is integer} \tag{24}$$

To prove the optimality conditions of the whole SC profit function, we firstly take the first and second partial derivatives of $\Pi_{sc/c}(p,Q,n)$ w.r.t. $n$ for fixed values of $p$ and $Q$ as follows:

$$\frac{\partial \Pi_{sc/c}}{\partial n} = \frac{(1-b)(\alpha-\beta p+\lambda\theta p)A_m}{(1-k^{1-b})Q^{1-b}n^2} - \frac{h_m(1-k)Q}{2}\left[1 - \frac{(1-b)(\alpha-\beta p+\lambda\theta p)(1-k)Q^b}{R(1-k^{1-b})}\right] \tag{25}$$

And

$$\frac{\partial^2 \Pi_{sc/c}}{\partial n^2} = -\frac{2(1-b)(\alpha-\beta p+\lambda\theta p)A_m}{(1-k^{1-b})Q^{1-b}n^3} < 0. \tag{26}$$

Clearly, $\Pi_{sc/c}(p,Q,n)$ is concave in $n$ for fixed values of $p$ and $Q$. Thus, like the procedure implemented in decentralized mode, the centralized model can be reduced to the calculation of a local optimal solution. Given $n$, taking the second derivative of $\Pi_{sc/c}$ w.r.t. $p$, we can find the SC profit function is concave on $p$ for a fixed value of $Q$, as shown in the following equation.

$$\frac{\partial^2 \Pi_{sc/c}}{\partial p^2} = -\frac{2(\beta-\lambda\theta)(1-\theta)(1-b)(1-k)}{(1-k^{1-b})}Q^b < 0 \tag{27}$$

Hence, we can get the optimal value of $p$, which is a function of $Q$, from $\partial\Pi_{m/d}/\partial n = 0$:

$$\frac{\partial \Pi_{sc/c}}{\partial p} = \frac{(1-b)(1-k)}{(1-k^{1-b})}Q^b\left[(1-\theta)\alpha + (\beta-\lambda\theta)\left(m - 2(1-\theta)P + \frac{\left(A_r+\frac{A_m}{n}\right)}{(1-k)Q} + \frac{h_m(2-n)(1-k)Q}{2R}\right)\right] = 0 \tag{28}$$

Which leads to

$$p^{**}(Q) = \frac{1}{2}\left[\frac{\alpha}{(\beta-\lambda\theta)} + \frac{1}{(1-\theta)}\left(m + \frac{\left(A_r+\frac{A_m}{n}\right)}{(1-k)Q} + \frac{h_m(2-n)(1-k)Q}{2R}\right)\right] \tag{29}$$

Substituting $p^{**}(Q)$ into Eq. (24), we can replace the SC profit function as a function of $Q$, as follows:

$$\Pi_{sc/c}(Q) = \frac{(\beta-\lambda\theta)(1-b)(1-k)}{4(1-k^{1-b})}\bigg[\left((1-\theta)\rho^2 + \frac{6}{(1-\theta)}\hat{A}\hat{H}\right)Q^b + \frac{3}{(1-\theta)}\left(\frac{\hat{A}^2}{(1-k)^2}Q^{b-2} + \hat{H}^2(1-k)^2Q^{b+2}\right) - 2\rho\left(\frac{\hat{A}}{(1-k)}Q^{b-1} + \hat{H}(1-k)Q^{b+1}\right)\bigg] - \left(\frac{(1-b)(1-k^{2-b})h_r}{(2-b)(1-k^{1-b})} + \frac{(1-k)(n-1)h_m}{2}\right)Q \quad 0 \leq \theta < \frac{\beta}{\lambda} < 1 \tag{30}$$

Where $\rho = \frac{\alpha}{(\beta-\lambda\theta)} - \frac{m}{(1-\theta)}$, $\hat{A} = A_r + \frac{A_m}{n}$ and $\hat{H} = \frac{h_m(2-n)}{2R}$. Meanwhile, simplifying Eq. (28), it can be illustrated as a convex set, as shown in Eq. (31):

$$\Pi_{sc/c}(Q) = \frac{(\beta-\lambda\theta)(1-b)(1-k)Q^b}{4(1-k^{1-b})}\bigg[(1-\theta)\rho^2 - 2\rho\left(\frac{\hat{A}}{(1-k)Q} + \hat{H}(1-k)Q\right) + \frac{3}{(1-\theta)}\left(\frac{\hat{A}}{(1-k)Q} + \hat{H}(1-k)Q\right)^2\bigg] - \left(\frac{(1-b)(1-k^{2-b})h_r}{(2-b)(1-k^{1-b})} + \frac{(1-k)(n-1)h_m}{2}\right)Q \tag{31}$$

Hence, from Eq. (32), the economic order quantity $Q^{**}$ can be obtained by setting $\partial \Pi_{sc/c}(Q)/\partial Q = 0$:

$$\frac{(\beta-\lambda\theta)(1-b)(1-k)}{4(1-k^{1-b})}\left[b\left((1-\theta)\rho^2 + \frac{6}{(1-\theta)}\hat{A}\hat{H}\right)Q^{b-1} + \frac{3}{(1-\theta)}\left((b+2)\hat{H}^2(1-k)^2Q^{b+1} - (2-b)\frac{\hat{A}^2}{(1-k)^2}Q^{b-3}\right) - 2\rho\left((b+1)\hat{H}(1-k)Q^b - (1-b)\frac{\hat{A}}{(1-k)}Q^{b-2}\right)\right] - \left(\frac{(1-b)(1-k^{2-b})h_r}{(2-b)(1-k^{1-b})} + \frac{(1-k)(n-1)h_m}{2}\right) = 0 \quad (32)$$

Meanwhile, we can apply a solution algorithm to find the optimal values of decision variables as follows.

**Step1:** Set $n = 1$, the minimum feasible value for $n$.

**Step2:** To achieve the maximum SC profit, $\Pi^*_{sc/c}(n, p^{**}(n), Q^{**}(n))$, calculate iteratively $p^{**}(n)$ and $Q^{**}(n)$ using Eqs. (29) and (32) until the solutions converge.

**Step3:** Set $n = n + 1$, and go to Step 2 in order to get new $\Pi^*_{sc/c}(n, p^{**}(n), Q^{**}(n))$.

**Step4:** If $\Pi^*_{sc/c}(n, p^{**}(n), Q^{**}(n)) > \Pi^*_{sc/c}((n-1), p^{**}(n-1), Q^{**}(n-1))$, then go to Step 3; otherwise, go to Step 5.

**Step5:** Set $\Pi^*_{sc/c}(n, p^{**}(n), Q^{**}(n)) = \Pi^*_{sc/c}((n-1), p^{**}(n-1), Q^{**}(n-1))$ where $n^{**}, p^{**}(n^{**}), Q^{**}(n^{**})$ represent a local optimum solution.

According to the Proposition 1, we prove that the solution derived from above iterative algorithm is, indeed, the global optimum solution.

**Proposition 1.** $\Pi^*_{sc/c}(n^{**}, p^{**}(n^{**}), Q^{**}(n^{**}))$ obtained by the iterative procedure outlined above is the global optimum solution.

**Proof.** It is found that the local optimum result calculated by the above mentioned procedure satisfies $\Pi^*_{sc/c}(n^{**}, p^{**}(n^{**}), Q^{**}(n^{**})) > \Pi^*_{sc/c}((n^{**}+1), p^{**}(n^{**}+1), Q^{**}(n^{**}+1))$, as shown in Fig. 3.

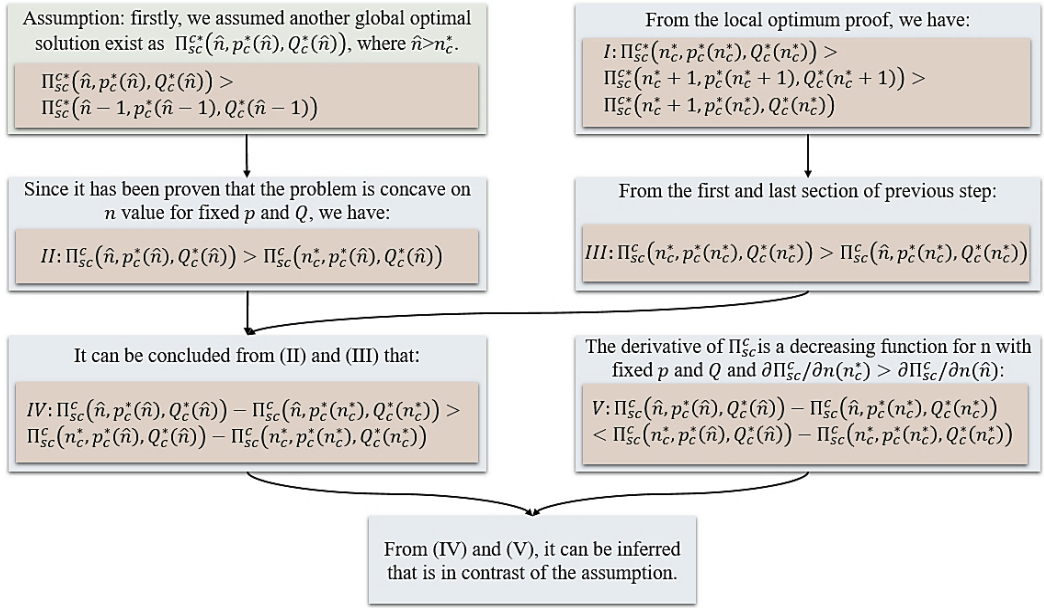

**Fig. 3.** Diagram of the global optimal solution proof.

Therefore, from Fig. 3, it can be found that the global optimum solution is the same $\Pi^*_{sc/c}(n^{**}, p^{**}(n^{**}), Q^{**}(n^{**}))$. It should be noted that the algorithm proposed above is applied for determining $n^{**}$.

Although decision making under the centralized system optimize and maximize the whole SC's profit, it may lead the retailer to face losses. In this vein, the retailer becomes reluctant to accept the optimal decisions adopted in the centralized system and prefer to operate independently based on the decisions he/she made under the decentralized system. Therefore, in order to incentivize the retailer to accept the global optimum solutions obtained in the centralized system, a motivating mechanism upon a combined revenue and cost-sharing (RCS) contract is developed in the next section to coordinate the SC members' decisions.

### 3.4. Coordinated system (RCS)

In this section, in order to achieve a win-win output for SC members, we propose and design a combined revenue and cost-sharing (RCS) contract to stimulate the retailer in adopting the centralized' global optimum decisions. To be more specific, the manufacturer decreases the wholesale price with a certain discount amount to make the retailer motivated to participate in this incentive scheme. On the other hand, the retailer is receptive to share $(1 - \mu)$ fraction of its revenue with the manufacturer only if the latter is satisfied to pay the former's holding cost in the same proportion. In this regard, to attain the perfect coordination decision-making system, we seek to determine the optimal values of new wholesale price and $\mu$. Therefore, the profit function of the retailer and the manufacturer under the proposed contract can be modelled as follows:

$$\Pi_{r/co}(p, Q, \mu) = \frac{1}{T_r}\left[(\mu p - v)(Q - kQ) - A_r - \mu h_r \int_0^{T_r} q(t)dt\right] \quad (33)$$

And

$$\Pi_{m/co}(n, \mu) = \frac{1}{T_r}\left\{(v - m - \theta p)(1 - k)Q - \frac{A_m}{n} - \frac{1}{2}h_m T_r\left[(n-1)\left(1 - \frac{n(1-k)Q}{RT}\right) + \left(\frac{n(1-k)Q}{RT}\right)\right](1 - k)Q + (1 - \mu)p(1 - k)Q - (1 - \mu)h_r \int_0^{T_r} q(t)dt\right\} \quad (34)$$

respectively. Substituting Eq. (5) into Eq. (33) and after simplifying, we have

$$\Pi_{r/co}(p, Q, \mu) = \frac{(1-b)(\alpha - \beta p + \lambda \theta p)}{(1-k^{1-b})}\left[(\mu p - v)(1-k)Q^b - A_r Q^{b-1}\right] - \mu \frac{(1-b)(1-k^{2-b})h_r}{(2-b)(1-k^{1-b})}Q \quad (35)$$

Taking the second derivative of $\Pi_{r/ci}$ w.r.t. $p$, we can be found the retailer profit function is concave on $p$ for a fixed value of $Q$, as shown in following.

$$\frac{\partial^2 \Pi_{r/co}}{\partial p^2} = -\frac{2(\beta - \lambda \theta)(1-b)(1-k)\mu}{(1-k^{1-b})}Q^b < 0 \quad (36)$$

Hence, the optimal retail price $p_{co}$ as a function of the order quantity $Q$ can be obtained at $\partial \Pi_{r/co}/\partial p = 0$.

$$p_{co}(Q_{co}) = \frac{1}{2}\left[\frac{\alpha}{\beta - \lambda \theta} + \frac{1}{\mu}\left(v + \frac{A_r}{(1-k)Q_{co}}\right)\right] \quad (37)$$

Applying the coordinating condition, where the optimal decision variables of the coordination system are equal to those of the centralized system ($p_{co} = p^{**}$, $Q_{co} = Q^{**}$, $n_{co} = n^{**}$),

The new discounted wholesale price, $v_{co}$, can be extracted from Eq. (37) as following

$$v_{co} = \frac{1}{(1-\theta)} \left[ \mu \left( m + \frac{A_r + \frac{A_m}{n^{**}}}{(1-k)Q^{**}} + \frac{h_m(2-n^{**})(1-k)Q^{**}}{2R} \right) - \frac{(1-\theta)A_r}{(1-k)Q^{**}} \right] \tag{38}$$

and consequently, the discount rate can be defined by $d = 1 - (v_{co}/v)$.

Substituting Eq. (38) into Eqs. (34) and (35), and after simplifying, we have

$$\Pi_{r/co}(p^{**}, Q^{**}, \mu) = \frac{(1-b)}{(1-k^{1-b})} \mu \eta \text{ (36)} \tag{39}$$

and

$$\Pi_{m/co}(n^{**}, \mu) = (1 - \theta - \mu)\eta - \left( \frac{(1-b)(1-k^{2-b})h_r}{(2-b)(1-k^{1-b})} + \frac{(1-k)(n^{**}-1)h_m}{2} \right) Q^{**} \tag{40}$$

, where

$$\eta = (\alpha - (\beta - \lambda\theta)p^{**})(1-k)Q^{**b} \left[ p^{**} - \frac{1}{1-\theta} \left( m + \frac{A_r + \frac{A_m}{n^{**}}}{(1-k)Q^{**}} + \frac{h_m(2-n^{**})(1-k)Q^{**}}{2R} \right) \right] - \frac{(1-k^{2-b})h_r}{(2-b)} Q^{**}. \tag{41}$$

3.4.1. The retailer's condition

The retailer will accept to participate in the incentive scheme proposed in the coordination system only if its profit under this system was not less than that under the decentralized system. Thus, the retailer's condition in order to agree with the contract is as follows:

$$\Pi_{r/co}(p^{**}, Q^{**}, \mu) \geq \Pi_{r/d}(p^*, Q^*) \tag{42}$$

By simplifying it, a lower limit on the variable decision $\mu$ will be obtained from the retailer's perspective as:

$$\underline{\mu} = \frac{1}{\eta} \left\{ (\alpha - (\beta - \lambda\theta)p^*)(1-k)Q^{*b} \left[ p^* - \left( v + \frac{A_r}{(1-k)Q^*} \right) \right] - \frac{(1-k^{2-b})h_r}{(2-b)} Q^* \right\} \tag{43}$$

3.4.2. The manufacturer's condition

Akin to the retailer's position, the manufacturer will accept to participate in the incentive scheme proposed in the coordination system only if its profit under this system was not less than that under the decentralized system. Thus, the manufacturer's condition in order to agree with the contract is as follows:

$$\Pi_{m/co}(n^{**}, \mu) \geq \Pi_{m/d}(n^*) \tag{44}$$

From Eq. (44), an upper limit on the variable decision $\mu$ will be obtained from the manufacturer's perspective as:

$$\bar{\mu} = 1 - \theta - \frac{1}{\eta}\left\{(\alpha - (\beta - \lambda\theta)p^*)(1-k)Q^{*b}\left[v - \left(\theta p^* + m + \frac{\frac{A_m}{n^*}}{(1-k)Q^*} + \frac{h_m(2-n^*)(1-k)Q^*}{2R}\right)\right] + \frac{(1-k^{2-b})h_r}{(2-b)}Q^{**} - \frac{h_m(1-k)(1-k^{1-b})h_m}{2(1-b)}\left((n^*-1)Q^* - (n^{**}-1)Q^{**}\right)\right\} \quad (45)$$

If $\bar{\mu} - \underline{\mu} \geq 0$, then supply chain coordination is feasible. Since $\vartheta \in [\underline{\mu}, \bar{\mu}]$, setting $\vartheta$ at any range of $[\underline{\mu}, \bar{\mu}]$ leads to an increase in the members' profits compared to the decentralized decision-making model, and the amount of each member's share from this additional profit depends on the bargaining power of each member against that of other members.

### 3.4.3. Profit allocation policy through bargaining

In this section, a profit allocation policy developed based on the SC members' bargaining power is presented. Indeed, we adopt the bargaining power as a criteria to adjust $\mu$ in revenue and cost-sharing contract. The bargaining of the retailer is considered $\xi$, where $0 < \xi < 1$, and therefore that of the manufacturer is equal to $1 - \xi$. By utilizing the proposed profit allocation policy, the extra profit resulted from the coordination system (denoted by $\Delta\Pi$) is divided between the SC members based on their bargaining power. Under this kind of decision-making system, the retailer's profit in the coordination decision-making system, $\Pi_{r/co}(p^{**}, Q^{**})$, is equal to its profit in the decentralized one $\Pi_{r/co}(p^{**}, Q^{**})$ plus a fraction of the extra profit, determined according to its bargaining power. Therefore, the retailer's profit under profit allocation policy is

$$\Pi_{r/co}(p^{**}, Q^{**}) = \Pi_{r/d}(p^*, Q^*) + \xi\Delta\Pi \quad (46)$$

As a result of this, the manufacturer seizes the remaining extra profit, $(1-\xi)\Delta\Pi$. As the procedure presented above for the retailer, the manufacturer's profit under coordination decision-making system with the profit allocation policy is

$$\Pi_{m/co}(n^{**}) = \Pi_{m/d}(n^*) + (1-\xi)\Delta\Pi \quad (47)$$

Where

$$\Delta\Pi = \Pi_{sc/co}(p^{**}, Q^{**}, n^{**}) - \Pi_{sc/d}(p^*, Q^*, n^*). \quad (48)$$

According to Eq. (46), substituting Eq. (7) into it and simplifying, the accurate factor of revenue and cost-sharing, $\mu_B$, can be calculated upon the bargaining power as follows:

$$\mu_B = \xi\bar{\mu} + (1-\xi)\underline{\mu} \quad (49)$$

, which this value is desirable to all members and allocates the extra profit fairly to them.

## 4. Numerical results and sensitive analysis

In this section, five numerical experiments as test problems are presented in order to evaluate the proposed model under different decision-making systems.

5.1. Numerical experiments and sensitively analysis

The Table. 2 presents data provided for five test problems.

Table 2. Five test problem parameters.

| Problem # | $\alpha$ | $(\beta, b, \lambda)$ | $\theta$ | $k$ | $R$ | $(v, m)$ | $(A_r, A_m)$ | $(h_r, h_m)$ | $\xi$ |
|---|---|---|---|---|---|---|---|---|---|
| 1 | 1200 | (8, 0.1, 9) | 0.15 | 0.6 | 1600 | (45, 10) | (250, 500) | (10, 5) | 0.4 |
| 2 | 900 | (12, 0.2, 15) | 0.2 | 0.4 | 6500 | (40, 25) | (100, 250) | (12, 7) | 0.5 |
| 3 | 1600 | (14, 0.2, 16) | 0.2 | 0.5 | 2100 | (70, 30) | (300, 450) | (15, 6) | 0.6 |
| 4 | 500 | (10, 0.3, 14) | 0.3 | 0.5 | 5000 | (50, 20) | (150, 300) | (9, 4) | 0.5 |
| 5 | 2500 | (15, 0.1, 18) | 0.15 | 0.6 | 8000 | (50, 15) | (200, 400) | (20, 10) | 0.4 |

Table. 3 shows the results run for the three proposed decision-making systems. In all experiments, the centralized decision-making system leads to a considerable increase in the profit for both the whole SC and the manufacturer in comparison with the decentralized system, while the retailer experiences a drop in its profit under the centralized system compared to the decentralized one. In this regard, the retailer is more likely to avoid adopting the decisions made in the centralized system and therefore, the manufacturer has no choice but to motivate the retailer to accept the optimal decisions of the centralized system by providing an incentive scheme proposed in the coordination system. As can be seen in Table. 3, a growth took place for the profits of the whole SC and SC members under the coordination decision-making system compared to the decentralized system. Therefore, the revenue and cost-sharing contract ensure that adopting the optimal decisions of the centralized system are beneficial for the SC to participate in such a decision-making system. Also, the extra profit resulted from the centralized system is divided according to the members' bargaining power. By considering the optimal decisions of the centralized system for the coordination system, the wholesale price is decreased with a discount rate, demonstrating the optimal performance of the revenue and cost-sharing contract in achieving the channel coordination.

Table 3. Results of the proposed models for five investigated test problems.

|  | Problem 1 | Problem 2 | Problem 3 | Problem 4 | Problem 5 |
|---|---|---|---|---|---|
| **Decentralized sys** | | | | | |
| $Q^*$ | 803.393 | 688.222 | 1205.16 | 552.893 | 930.268 |
| $p^*$ | 113.11 | 70.12 | 109.32 | 68.37 | 126.9 |
| $n^*$ | $\lceil 1.88 \rceil = 2$ | $\lceil 0.66 \rceil = 1$ | $\lceil 1.73 \rceil = 2$ | $\lceil 1.18 \rceil = 1$ | $\lceil 1.16 \rceil = 1$ |
| $\Pi_{r/d}$ | 51079.8 | 21716.92 | 49766.5 | 7476.15 | 123908 |
| $\Pi_{m/d}$ | 13930.7 | 136.05 | 27118.5 | 5194.17 | 26634.6 |
| $\Pi_{sc/d}$ | 65010.6 | 21852.97 | 76885 | 12670.32 | 150542.6 |
| **Centralized sys** | | | | | |
| $Q^{**}$ | 1007.78 | 754.621 | 2457.64 | 1196.29 | 1229.03 |
| $p^{**}$ | 96.83 | 66.26 | 89.73 | 58.09 | 111.34 |
| $n^{**}$ | 2 | 1 | 5 | 1 | 1 |
| $\Pi_{r/c}$ | 47497.7 | 21232.21 | 27617.3 | 1773.25 | 117430.2 |
| $\Pi_{m/c}$ | 20527.6 | 1005.09 | 62238.5 | 14725.35 | 38637.8 |
| $\Pi_{sc/c}$ | 68025.3 | 22237.3 | 89855.8 | 16498.6 | 156068 |
| **Coordinated sys** | | | | | |
| $Q^{**}$ | 1007.78 | 754.621 | 2457.64 | 1196.29 | 1229.03 |
| $p^{**}$ | 96.83 | 66.26 | 89.73 | 58.09 | 111.34 |
| $n^{**}$ | 2 | 1 | 5 | 1 | 1 |
| $\mu$ | 0.618 | 0.74 | 0.361 | 0.276 | 0.662 |
| $\bar{\mu}$ | 0.654 | 0.752 | 0.455 | 0.418 | 0.691 |
| $\mu_B$ | 0.632 | 0.746 | 0.417 | 0.347 | 0.673 |
| $v_{co}$ | 7.73 | 24.03 | 12.81 | 10.15 | 12.67 |
| $d(\%)$ | 82.82 | 39.92 | 81.7 | 79.7 | 74.66 |
| $\Pi_{r/co}$ | 52225.4 | 21903.4 | 57384.1 | 9365.76 | 125899 |
| $\Pi_{m/co}$ | 15799.9 | 333.9 | 32471.7 | 7132.84 | 30168.3 |
| $\Pi_{sc/co}$ | 68025.3 | 22237.37 | 89855.8 | 16498.6 | 156068 |
| $Savings^i$ (%) | | | | | |
| Retailer | 2.24 | 0.85 | 15.3 | 25.27 | 1.6 |
| Manufacturer | 13.41 | 145.42 | 19.74 | 37.32 | 13.26 |
| SC | 4.63 | 1.75 | 16.87 | 30.21 | 3.67 |

$^i Savings = ((\Pi_{co} - \Pi_d)/\Pi_d) \times 100\%$.

The developed models are analyzed by the key parameter $\theta$, the portion of retail price considered for the social donation activity, in order to evaluate them under varied decision-making systems. The reason of selecting this parameter for the examination is that it has the highest impact on variable decisions and the profit of the whole SC and the members due to the significance of social activity in seizing a higher market share. The data of test problem (1) are applied for the sensitively analysis.

Figure 4 shows the effect of changing $\theta$ on the retail price. It is clear that the more the portion of retail price donated to social activity $\theta$, the higher the value of retail price under all three decision-making systems, as shown in Fig. 4. In fact, the SC prefer to increase the retail price when it is supposed to allot a higher portion of sales to social activity. Though, the growth in the retail price under the decentralized system is more than that under the centralized/coordination systems, which increases the market demand share resulting in the higher profit under centralized/coordination system.

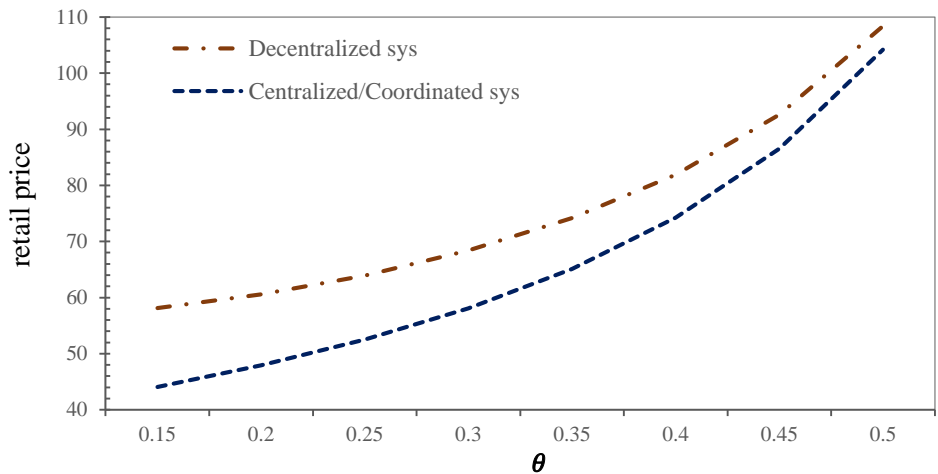

**Fig 4**. Impact of donated portion of retail price on optimal retail price

According to Fig. 5, the retailer's order quantity increases when $\theta$ increases because the market demand will increase by increasing $\theta$ through a cause-related marketing (CRM) campaign, stimulating the retailer to order more in order to satisfy rising demand and prevent from facing shortage. It can be seen that the order quantity in the centralized/coordination system is less than that in the decentralized system when the value of $\theta$ becomes more than 0.4. It can be rooted in this fact that a single decision maker in the centralized/coordination system avoid ordering fearlessly due to being well aware of the SC members' conditions in the market; however, an independent retailer under the decentralized decision-making system order promptly while facing much more demand regardless of whether or not it is profitable for the whole SC.

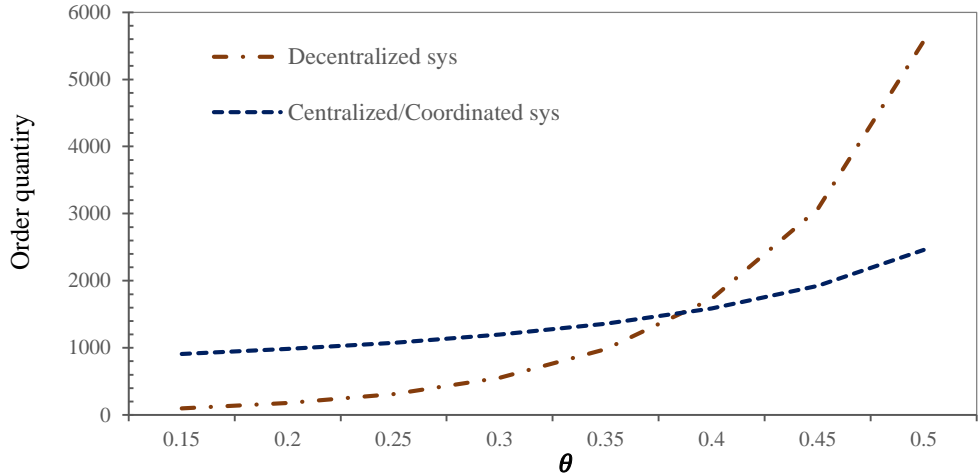

**Fig 5.** Impact of donated portion of retail price on economic order quantity

Fig. 6 demonstrates the changes in the bounds of the incentive RCS contract with the changes of $\theta$. As we mentioned in 3.4.2., the SC would be coordinated when $\overline{\mu}$ is always higher than $\underline{\mu}$. Clearly, the performance of the incentive contract is possible for each value of $\theta$. Indeed, it is clear that the RCS contract is efficient in coordinating the SC under each value of $\theta$. Moreover, the achievement to a coordination system is easier when the value of $\theta$ is lower because considering a higher proportion of retail price for SD activity is likely to create loss for the manufacturer (see Fig. 8).

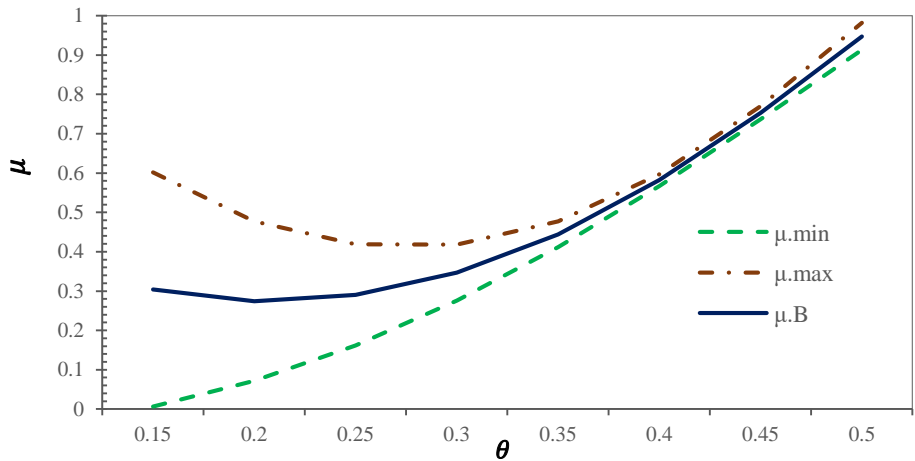

**Fig 6.** $\underline{\mu}$ versus $\overline{\mu}$ for alternative values of $\theta$.

Figures (7) and (8) demonstrate the impacts of various $\theta$ on the expected profits of the retailer and the manufacturer. As can be seen, the profit of the retailer in the centralized system decreased compared in the decentralized system while that of the manufacturer experienced a rise in the centralized system. As mentioned, the retailer avoids adopting the global optimal decisions made under the centralized system while the manufacturer is willing to participate in this system. As the profit of both manufacturer and whole SC increased in the centralized system compared to the decentralized, the manufacturer makes an effort to persuade the retailer to accept the global optimal decisions by offering the incentive contract in the coordination system. As can be seen in the both figures, the figures of the coordinated system are always upper than those of the decentralized system for both members; therefore, accepting the global optimum decisions with a balancing contract under the coordination system is beneficial for them. Moreover, by increasing the portion of retail price $\theta$ allocated to SD activity, the retailer's profit will increase while the manufacturer's will decreases. It is due to the fact that the demand imposed on the retailer will rise while a growth take place on the portion of retail price $\theta$ allocated to SD activity since customers concerned about social issues prefer to buy goods with a higher portion allotted to social activities (see Fig. 7). However, because the manufacturer accepted to be responsible for paying costs related to social issues, he/she experience a declining trend in his profit due to allocating more budget to this issue by increasing $\theta$ (see Fig. 8). According to Fig. 8, operating in the SD activity is feasible for the manufacturer only when $\theta <$ 0.4, otherwise he/she suffers losses because there always is a budget limit for participating in any activity.

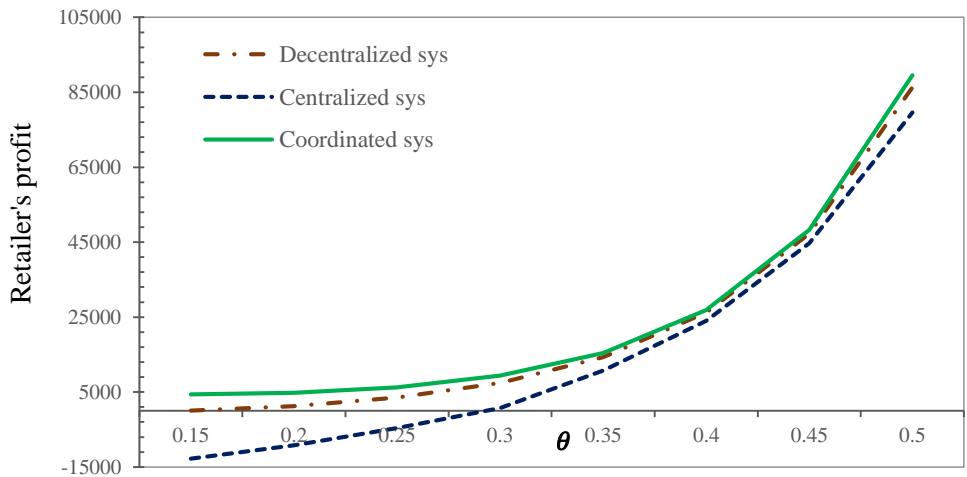

**Fig 7**. Impact of donated portion of retail price on the retailer's profit

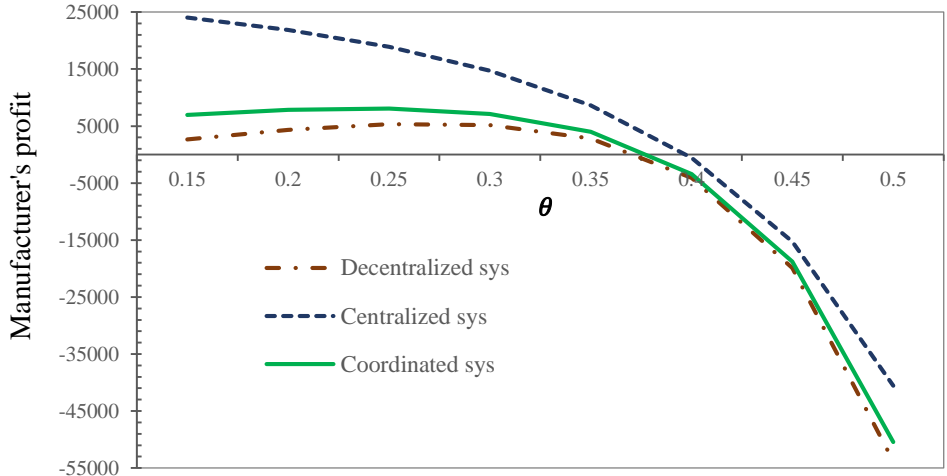

**Fig 8**. Impact of donated portion of retail price on the manufacturer's profit

According to Fig. 9, the whole SC profit increases with a rise in $\theta$, and it is clear from the whole SC function which is a growing function of $\theta$ under various systems. Moreover, the line graph of the coordinated system is always upper than those of the decentralized system, meaning that operating under the coordination system is achievable for all amount of $\theta$.

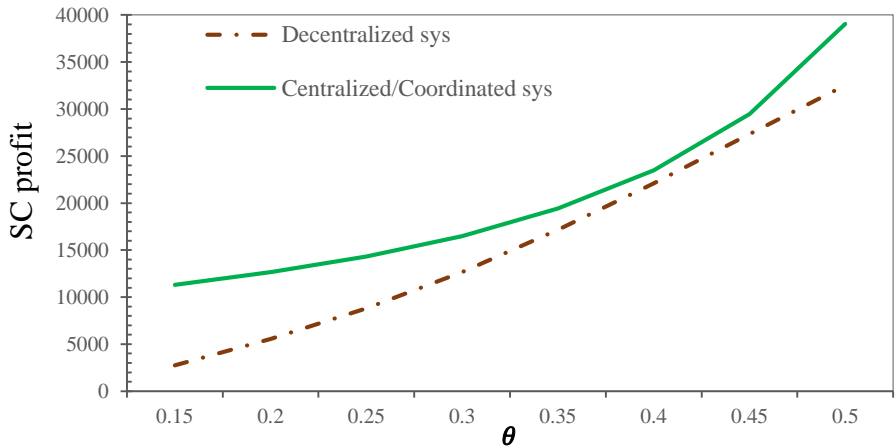

**Fig 9**. Impact of donated portion of retail price on the SC profit

## 5. Blocked decision-making on social donation (SD) activity

In the previous section, simultaneous decision making was developed on pricing, replenishment and SD decisions with the performance of incentive RCS mechanism. However, pricing, inventory and social activity decisions are likely to be made individually by different entities of each organization. In this regard, in this section the model is analyzed under a circumstance in which the decisions on pricing and inventory are adopted independently from decision making on SD activity. This is to emphasize the profitability and efficiency of the previous proposed models and shows how the performance of the supply chain with and without taking into account the SD decision.

By neglecting the impact of SD activity on the problem, all coefficients related to SD become equal to zero. Therefore, the previous demand function, $(\alpha - \beta p + \lambda \theta p)q(t)^b$, will be converted into $(\alpha - \beta p)q(t)^b$, resulting in $T_r = (1 - k^{1-b})Q^{1-b}/(1 - b)(\alpha - \beta p)$. The decentralized and centralized systems first are modelled and then, the revenue and cost-sharing (RCS) mechanism will be implemented to achieve a coordination system. Finally, the results of the proposed model without considering social activity are analyzed and compared with the developed model in which the SD activity was considered.

5.1. Decentralized system

Under the decentralized system, the retailer's profit function is modelled as follows:

$$\Pi_r(p.Q) = \frac{(1-b)(\alpha - \beta p)}{(1 - k^{1-b})}[(p-v)(1-k)Q^b - A_r Q^{b-1}] - \frac{(1-b)(1-k^{2-b})h_r}{(2-b)(1-k^{1-b})}Q \qquad p < \frac{\alpha}{\beta} \qquad (50)$$

Similar to the process of proving the concavity in section 3.2.1, the optimal value of retailing price $p^*$ is calculated as follows:

$$p^*(Q) = \frac{1}{2}\left[\frac{\alpha}{\beta} + v + \frac{A_r}{(1-k)Q}\right] \qquad (51)$$

Substituting Eq. (60) into Eq. (50), the retail's profit function is obtained as a function of $Q$ as follows,

$$\Pi_r(Q) = \frac{\beta(1-b)}{4(1-k^{1-b})}\left[\left[\left(\frac{\alpha}{\beta} - v\right)\sqrt{(1-k)Q^b}\right) - \frac{A_r}{\sqrt{(1-k)Q^{2-b}}}\right]^2 - \frac{4(1-k^{2-b})}{\beta(2-b)}h_r Q\right] \qquad (52)$$

Accordingly, the economic order quantity $Q^*$ can be generated.

Furthermore, in this model, the profit function of the manufacturer is designed as following,

$$\Pi_m(n) = \frac{(1-b)(\alpha-\beta p)}{(1-k^{1-b})}\left[(v-m)(1-k)Q^b - \frac{A_m}{n}Q^{b-1}\right] - \frac{h_m(1-k)Q}{2}\left[(n-1) + \frac{(1-b)(\alpha-\beta p)(2-n)(1-k)Q^b}{R(1-k^{1-b})}\right] \qquad (53)$$

As this profit function is a concave function over $n$, the proof of which is similar to what was done in the section 3.2.2, $n^*$ is:

$$n^*_{Desimal} = \sqrt{\frac{2RA_m(1-b)(\alpha-\beta p^*)Q^{*b}}{h_m(1-k)Q^{*2}\left(R(1-k^{1-b}) - (1-b)(\alpha-\beta p^*)(1-k)Q^{*b}\right)}} \qquad (54)$$

In next section, decisions in the blocked decision-making model are made based upon the perspective of the entire SC in order to maximize the whole SC profit.

5.2. Centralized system

Under this decision-making system, the profit function of the entire SC is modelled as:

$$maximize\ \Pi_{sc}(p.Q.n) = \frac{(1-b)(\alpha-\beta p)}{(1-k^{1-b})}\left[(p-m)(1-k)Q^b - \left(A_r + \frac{A_m}{n}\right)Q^{b-1}\right] - \frac{(1-b)(1-k^{2-b})h_r}{(2-b)(1-k^{1-b})}Q - \frac{h_m(1-k)Q}{2}\left[(n-1) + \frac{(1-b)(\alpha-\beta p)(2-n)(1-k)Q^b}{R(1-k^{1-b})}\right] \quad (55)$$

To calculate the optimal values for the decision variables of the whole SC and prove the concavity, we operate according to the same approaches we adopted in the centralized system of the previous model, i.e., Proposition 1 and Solution algorithm in section 3.3. In this vein, the optimal value of retail price $p$ is obtained as a function of $Q$ with given $n$.

$$p^{**}(Q) = \frac{1}{2}\left[\frac{\alpha}{\beta} + m + \frac{\left(A_r + \frac{A_m}{n}\right)}{(1-k)Q} + \frac{h_m(2-n)(1-k)Q}{2R}\right] \quad (56)$$

Substituting Eq. (56) into Eq. (55), the whole SC function is obtained as a function of $Q$ as follows,

$$\Pi_{sc}(Q) = \frac{\beta(1-b)(1-k)Q^b}{4(1-k^{1-b})}\left[\left(\phi - \frac{\widehat{A}}{(1-k)Q}\right)^2 + (\phi - \widehat{H}(1-k)Q)^2 + (2\widehat{A}\widehat{H} - \phi^2)\right] - (51)\left(\frac{(1-b)(1-k^{2-b})h_r}{(2-b)(1-k^{1-b})} + \frac{\beta(1-k)(n-1)h_m}{4}\right)Q \quad (57)$$

Where $\phi = \frac{\alpha}{\beta} - m$, $\widehat{A} = A_r + \frac{A_m}{n}$ and $\widehat{H} = \frac{h_m(2-n)}{2R}$. Accordingly, it can be inferred that Eq. (51) is a convex set, as a result of which the economic order quantity $Q^{**}$ can be calculated through $\partial \Pi_{sc}(Q)/\partial Q = 0$.

Similar to the model with SD activity, decision making under the centralized system is likely to be uneconomical; therefore, an incentive mechanism would be efficient to persuade him to decide based on the centralized solutions. Hence, in the next section, the RCS contract is applied in the blocked decision-making model.

5.3. Coordination system: the RCS mechanism

In this section, the manufacturer seeks to motivate the retailer to adopt the decisions made in the centralized system. Similar to the conditions of the SC members in Eqs. (42) and (44), the lower and upper bounds of the factor $\mu$ in the RCS mechanism are obtained as follows, respectively.

$$\underline{\mu} = \frac{1}{\delta}\left\{(\alpha - \beta p^*)(1-k)Q^{*b}\left[p^* - \left(v + \frac{A_r}{(1-k)Q^*}\right)\right] - \frac{(1-k^{2-b})h_r}{(2-b)}Q^*\right\} \quad (58)$$

And

$$\bar{\mu} = 1 - \frac{1}{\delta}\left\{(\alpha - \beta p^*)(1-k)Q^{*b}\left[v - \left(p^* + m + \frac{\frac{A_m}{n^*}}{(1-k)Q^*} + \frac{h_m(2-n^*)(1-k)Q^*}{2R}\right)\right]\right.$$
$$+ \frac{(1-k^{2-b})h_r}{(2-b)}Q^{**} \tag{59}$$
$$\left. - \frac{h_m(1-k)(1-k^{1-b})h_m}{2(1-b)}\left((n^*-1)Q^* - (n^{**}-1)Q^{**}\right)\right\}$$

The coordination is feasible for each value of the factor belonging to the interval $[\underline{\mu}, \bar{\mu}]$. The exact value of $\mu$ is obtained based on the SC members' bargaining power as:

$$\mu_B = \frac{1}{\delta}\left\{\xi\Delta\Pi + (\alpha - \beta p^*)(1-k)Q^{*b}\left[p^* - \left(v + \frac{A_r}{(1-k)Q^*}\right)\right] - \frac{(1-k^{2-b})h_r}{(2-b)}Q^*\right\} \tag{60}$$

Where

$$\delta = (\alpha - \beta p^{**})(1-k)Q^{**b}\left[p^{**} - \left(m + \frac{A_r + \frac{A_m}{n^{**}}}{(1-k)Q^{**}} + \frac{h_m(2-n^{**})(1-k)Q^{**}}{2R}\right)\right]$$
$$- \frac{(1-k^{2-b})h_r}{(2-b)}Q^{**} \tag{61}$$

By applying $\mu_B$, the surplus profit of the coordination system is divided between the members according to their bargaining power.

In the next section, the results of the blocked decision-making model are analyzed and compared with those of the model with the SD factor.

5.4. Discussion based upon numerical results

To analyze the performance of the blocked decision making on SD activity, we utilize the values and parameters used for the test problem 1, outputs of which are generated in Table 4.

In the following table, the whole SC's profit under the centralized system increased remarkably in comparison to that under the decentralized one, though the retailer's profit decreased by switching from the decentralized system to the centralized system. By implementing the RCS incentive contract to achieve the coordination system, the whole SC's profit remains on the same optimum level (the profit obtained from the centralized system), and both SC members benefit from participating in this incentive mechanism compared to operating under the traditional decision-making structure (decentralized system).

Table 4. Outcomes of running the blocked decision-making model for the parameters of test problem 1.

| system | $Q$ | $p$ | $n$ | $\Pi_r$ | $\Pi_m$ | $\Pi_{sc}$ |
|---|---|---|---|---|---|---|
| decentralized | 601.8 | 98.01 | [2.26] =2 | 35238.3 | 25564.5 | 60802.8 |
| centralized | 991.43 | 80.21 | 2 | 29680.4 | 36375.3 | 66055.6 |
| coordinated | 991.43 | 80.21 | 2 | 36718.8 | 29336.8 | 66055.6 |

Comparing the results of running the blocked decision-making model with those of the model with SD decisions, it can be inferred that applying the latter model results in more profits for the whole SC and its members than the former model. To be more precise, simultaneous decision making of SD, pricing and replenishment decisions leads to up to a 3% grow in the whole SC's profit, compared to the case in which the correlation between the social, pricing and inventory decisions is overlooked. From another perspective, when the retailer decides on pricing and inventory decisions regardless of considering the SD activity, the values of the retail price and order quantity would be different from those in the model with the social factor. For instance, in the main proposed model, the retail price is equal to 96.83 while this value decreases to 80.21 by adopting blocked decision-making model. It means that the price of the goods contributing to the social activity (making their customers informed by sticking some labels on the products) is higher than that of goods without any social activity. Surprisingly, despite an increase in the retail price, the market demand increases due to the motivation of customers to buy such products, resulting in more economic order quantity by the retailer.

In conclude, since pricing, inventory and SD decisions are very influential on the whole SC's sales volume, joint and simultaneous decision-making of such important factors is more likely to make much more profit for the whole SC and its members. Furthermore, applying coordination models along with taking into account the correlation between such issues leads to much more productive improvements in SC's performance.

## 7. Conclusion

In this article, simultaneous coordination of pricing, replenishment and social activity decisions in a two-echelon SC, including one manufacturer and one retailer, with market demand dependent on the retail price, the retailer's inventory level and social donation (SD) activity was investigated. While the retailer decided on the price and order quantity, the manufacturer aimed to engage in the SD activity through consumer social awareness (CSA) by labelling the products to declare that a proportion of the sales is allocated to the SD activity. By adopting such strategies, SC intended to seize more market share and profit, resulting from facing much higher market demand. The problem was modelled under three different decision-making systems. Firstly, the decentralized model, in which SC members were willing to maximize their profit by independently decision making regardless of the profits of the whole SC and its members, was proposed. Then, we presented a model under the centralized decision-making system, where optimal solutions from the viewpoint of whole SC were obtained. Although the centralized system results in the most expected profit (optimal profit) for the whole SC, it led the retailer to experience a decrease in his profit compared to his decentralized profit. Under such a circumstance, to convince the retailer to accept the optimal solutions obtained in the centralized system, we proposed a motivating scheme based on the RCS contract. Indeed, this mechanism incentivize both SC members to participate in the coordination system. In order to fairly divide the surplus profit gained from the centralized system between the members, a profit-sharing strategy based on the bargaining power of the SC members was presented, which is similar to what happens in the real-world business. The outputs of running models illustrate that SD activities can be an applicable strategy in seizing much higher market demand, and the performance of the SC with SD activity became more efficient than that of blocked decision-making model, which neglected the impact of SD activity. Moreover, the coordination of such decisions led to benefits for the whole SC and its members Although participating in SD activity was costly for the manufacturer (see Fig. 8), the whole SC and its members' profit increased by operating under the coordinated decision-making system.

Akin to any other research, there are some operational facets for this article which can be considered as contributions and extensions for the model in future studies. This model can be developed for the cases with other types of SC network, including multi-echelon SC and multiple players/members at each echelon under different game structures. Furthermore, the demand function can have different distribution functions and be dependent on other topics by adjusting into real world data, such as promotional efforts, credit time and quality. The RCS contract was applied as final coordinator, while it can be improper for some real cases; therefore, using another contracts, such as credit, sale rebate, buy back, and so on, may be suitable for other cases in the real-world business. Finally, while this study considered social aspect as a sustainability perspective along with the economic aspect, environmental aspect as another sustainability option should also be adopted in such a model.

# References

This section intentionally left blank